\documentclass[MSNbibl,number,citesort,seceqn,dvips]{arxbj}
\usepackage{mathrsfs}

\aid{0}
\volume{20}
\issue{1}
\pubyear{2014}
\firstpage{377}
\lastpage{393}
\doi{10.3150/12-BEJ490} 

\makeatletter

\newcommand{\RMO}{\mathrm{O}}

\newcommand{\RMe}{\mathrm{e}}

\newcommand{\cal}{\mathcal}

\newtheorem{theorem}{Theorem}
\newtheorem{lem}{Lemma}

\def\P{{\mathbb P}}
\def\E{{\mathbb E}}
\def\I{{\mathbb I}}
\def\R{{\mathbb R}}

\newcommand{\cW}{{\cal W}}
\newcommand{\cZ}{{\cal Z}}

\makeatother

\begin{document}
\begin{frontmatter}

\title{Small value probabilities for supercritical branching
processes with immigration}
\runtitle{Small value probabilities}

\begin{aug}
\author[1,2]{\fnms{Weijuan} \snm{Chu}\thanksref{1,2}\ead[label=e1]{chuwj@nju.edu.cn}},
\author[3]{\fnms{Wenbo V.} \snm{Li}\thanksref{3}\ead[label=e2]{wli@math.udel.edu}} \and
\author[1,4]{\fnms{Yan-Xia} \snm{Ren}\corref{}\thanksref{1,4}\ead[label=e3]{yxren@math.pku.edu.cn}}
\runauthor{W. Chu, W.V. Li and Y.-X. Ren} 
\address[1]{LMAM School of Mathematical Sciences, Peking University,
Beijing 100871, P.R. China}
\address[2]{Department of Mathematics, Nanjing University, Nanjing
210093, P.R. China.\\ \printead{e1}}
\address[3]{Department of Mathematical Sciences,
University of Delaware, Newark, DE 19716, USA.\\ \printead{e2}}
\address[4]{Center for Statistical Science, Peking University, Beijing
100871, P.R. China.\\ \printead{e3}}
\end{aug}

\received{\smonth{3} \syear{2011}}
\revised{\smonth{7} \syear{2012}}

%
\begin{abstract}
We consider a supercritical Galton--Watson branching process with
immigration. It is well known that under suitable conditions on the
offspring and immigration distributions, there is a finite, strictly
positive limit $\cW$ for the normalized population size. Small value
probabilities for $\cW$ are obtained. Precise effects of the balance
between offspring and immigration distributions are characterized.
\end{abstract}

%
\begin{keyword}
\kwd{immigration}
\kwd{small value probability}
\kwd{supercritical Galton--Watson branching process}
\end{keyword}

\end{frontmatter}

\section{Introduction and main results}

Small value probability for a \textit{positive} random variable $V$
studies the rate of decay of the so called left tail probability $\P(V
\le
\varepsilon) $ as $\varepsilon\to0^+$. When $V$ is the norm of a
random element in a Banach space, one is dealing with small ball
probability, see
\cite{LS01} for a survey of Gaussian measures. When $V$ is the maximum
of a continuous random process starting at zero, one is
estimating lower tail probability which is closely related to studies
of boundary crossing probabilities or the first
exit time associated with a general domain, see \cite{L03} and \cite
{LS04} for Gaussian processes.
A~comprehensive study of small value probability is emerging and
available in various talks and lecture notes in
\cite{L12+}, see also the literature compilation \cite{sdbib}.

In this paper, we further study the most natural aspect of the
branching tree approach originated in \cite{MO08} on the martingale
limit of a supercritical Galton--Watson process. The problem has been
solved initially in \cite{D71a,D71b}, see Theorem \ref{THM}.
The main goal is developing additional
tools to treat small value probabilities for the martingale limit of
a supercritical Galton--Watson process with immigration. The
interplay between the offspring and the immigration distribution can
be seen clearly from our main result Theorem \ref{thm1}. We next provide a
more detailed and precise discussion by introducing additional
notations, surveying relevant results and stating our results.

Let $(Z_n,n\geq0)$ be a supercritical Galton--Watson branching
process with $Z_0=1$, offspring distribution $p_k=\P(X=k), k\geq0$,
and mean $m= \E X\in(1,\infty)$.
To avoid non-branching case, we
suppose $p_k<1$ for all $k$ throughout this paper. Under the natural
condition $\E[X \log^{+} X]<\infty$, the positive martingale $Z_n
m^{-n}$ converges to a non-trivial random variable $W< \infty$ in the
sense (see Kesten and Stigum \cite{KS66})
\[
Z_n m^{-n} \longrightarrow W \qquad\mbox{a.s. and }
L^1 \mbox{ as } n \to\infty.
\]
Here and throughout this paper, $\log^+ x=\log\max(x, 1)\ge0$. The
distribution of the limit $W$ is of great interests in various
applications. However, except for some very special cases, the explicit
distribution of $W$ is not available, see, for example, Harris
\cite{H48}, Hambly \cite{H95} and Williams \cite{W08}, Section 0.9. In
general, it is known that $W$ has a continuous positive density on
$(0,\infty)$ satisfying a Lipschitz condition, see Athreya and Ney
\cite{AN72}, Chapter II, page 84, Lemma 2. However, it is not clear what type
of densities can arise in this way. This lack of complete information
on the distribution of $W$ prompts a search for asymptotic information
such as the behavior of the left tail, or the small value probabilities
of $W$ and its density.

In \cite{D71b}, the following results were given with
assumption $p_0=0$ which holds
without loss of generality
after the standard Harris--Sevastyanov transformation, see \cite{H48},
page 478, Theorem~3.2, or \cite{B88}, page 216.
Here and throughout this paper, we use $g_1(x)\asymp g_2(x)$ as $x\to
0^+ (\infty)$ to represent $c\le g_1(x)/g_2(x)\le C$ as $x\to0^+
(\infty)$ for two constants $C>c>0$ and $g_1(x)\sim g_2(x)$ as $x\to
0^+ (\infty)$ to represent $g_1(x)/g_2(x)\to1$ as $x\to0^+
(\infty)$.
%
\begin{theorem}[(Dubuc \cite{D71b})]\label{THM}
\textup{(a)} If $p_1>0$, then
\[
\P(W\leq\varepsilon)\asymp\varepsilon^{|{\log p_1}|/\log m} \qquad\mbox{as }
\varepsilon\to
0^+.
\]

\textup{(b)} If $p_1=0$, then
\[
-\log\P(W\leq\varepsilon)\asymp\varepsilon^{-\beta/(1-\beta)} \qquad\mbox{as }
\varepsilon
\to0^+
\]
with $\beta:= \log\gamma/ \log m$ and $\gamma:=\inf\{n\dvtx  p_n >0\}
\geq2$.
\end{theorem}

Note that the rough asymptotic $\asymp$ in Theorem \ref{THM} cannot be
improved into more precise asymptotic $\sim$ and the oscillation is
very small. This is the so called near-constancy phenomenon that were
described and studied theoretically or numerically in
\cite{D82,B88,BP88} and \cite{BB91}. In fact, it is still an open
conjecture that the Laplace transform of $W$ being non-oscillating near
$\infty$ (and hence the small value probability of $W$ being
non-oscillating near $0^+$) is only specific to the case $p_1>0$ in
\cite{KM68a}, page 127. General estimates, near-constancy phenomena,
specific examples, and various implications have been studied to
various degree of accuracy in Harris \cite{H48}, Karlin and McGregor
\cite{KM68a,KM68b}, Dubuc \cite{D71a,D71b} and \cite{D82},
Barlow and Perkins \cite{BP88}, Goldstein \cite{G87}, Kusuoka \cite
{K87}, Bingham~\cite{B88}, Biggins and Bingham \cite{BB91} and
\cite{BB93}, Biggins and Nadarajah \cite{BN93}, Fleischman and Wachtel
\cite{FW07} and \cite{FW09}. Recently, Berestycki, Gantert, M\"orters
and Sidorova \cite{BGMS12} studied limit behaviors of the
Galton--Watson tree conditioned on $W<\varepsilon$ as
$\varepsilon\downarrow0$.

In the present paper, we consider the supercritical branching
process with immigration denoted by $(\cZ_n,n\geq0)$, and follow
the definition in \cite{AN72}, Chapter VI, Section 7.1, page 263. To be
more precise, we have
\[
\cZ_0=Y_0,\qquad \cZ_{n+1}=X_1^n+X_2^n+
\cdots+X_{\cZ
_{n}}^n+Y_{n+1},\qquad n\ge0,
\]
where $X_1^n,X_2^n,\ldots$ are i.i.d.
with the same offspring distribution,
$Y_0,Y_1,\ldots$ are i.i.d. with the same immigration distribution
$\{q_k,k\geq0\}$, and $X$'s and $Y$'s are independent.
Recall that the offspring number $X$ has distribution $p_k=\P(X=k),
k\ge0$ and mean $m=\E X$.
Suppose $Y$ has distribution $\{q_k, k\ge0\}$.
We use $f(s)$ and $h(s)$ to denote the generating function of $X$ and
$Y$, respectively, that is,
%
\begin{equation}
\label{fh}f(s)=\E s^X=\sum^\infty_{k=0}p_k
s^k \quad\mbox{and}\quad h(s)=\E s^Y=\sum
^\infty_{k=0}q_k s^k,\qquad 0<s<1.
\end{equation}
It is a classical result, see Seneta \cite{S70}, for example, that
%
\begin{equation}
\label{limitW} \lim_{n\to\infty}\cZ_n/m^n=\cW
\end{equation}
exists and is finite a.s. if and
only if
%
\begin{equation}
\label{XYlog} \E\log^+ Y < \infty\quad\mbox{and}\quad \E\bigl(X \log^{+} X
\bigr)<\infty.
\end{equation}

Our main result of this paper is the following small value
probabilities for $\cW$, which can be expressed as weighted summation
of an infinite independent sequence of
$W$'s, see~(\ref{identity}).
%
\begin{theorem}\label{thm1}
Assume the condition (\ref{XYlog}) holds.

\textup{(a)} If $p_0=0$ and $0<q_0<1$, then
%
\begin{equation}
\label{3}\P(\cW\leq\varepsilon)\asymp\varepsilon^{|{\log q_0}|/\log m}
\qquad\mbox{as }
\varepsilon\to0^+.
\end{equation}

\textup{(b)} If $p_0=0$, $q_0=0$ and $p_1>0$, then
%
\begin{equation}
\label{4} \log\P(\cW\leq\varepsilon)\sim-\frac{K |{\log p_1}| }{2 (\log
m)^2}\cdot|{\log
\varepsilon}|^2 \qquad\mbox{as } \varepsilon\to0^+
\end{equation}
with $K=\inf\{n\dvtx q_n>0\}$.

\textup{(c)} If $p_0=0$, $q_0=0$ and $p_1=0$, then
\[
\log\P(\cW\leq\varepsilon)\asymp-\varepsilon^{-\beta/(1-\beta
)} \qquad\mbox{as }
\varepsilon\to0^+
\]
with $\beta$ being defined as in Theorem \ref{THM}\textup{(b)}.

\textup{(d)} If $p_0>0$, then
%
\begin{equation}
\P(\cW\leq\varepsilon) \asymp\varepsilon^{|{\log h(\rho)}|/\log
m} \qquad\mbox{as } \varepsilon
\to0^+,
\end{equation}
where $\rho$ is the solution of $f(s)=s$ between $(0,1)$, $f$ and $h$
are defined in (\ref{fh}).
\end{theorem}

Note that there are additional phase transitions appearing in the case
with immigration, in particular between the case
where the immigration distribution has a positive mass at 0 and where
there is no mass at 0.
In the $p_0>0$ case, the extinction probability of the branching
process $(Z_n, n\ge0)$ (without immigration) is strictly positive,
and plays the dominating role in the small value probability of $\cW$.
Our approach is outlined in
Section \ref{lemma} and detailed proof of Theorem \ref{thm1} is give in
Sections \ref{lower}, \ref{upper} and \ref{p0}.

\section{Our approach}\label{lemma}

Our proof of Theorem \ref{thm1} is
based on Dubuc's result stated in Theorem \ref{THM}. In \cite{D71b},
an integral composition transform is used together with some
non-trivial complex analysis, which is powerful but inflexible and
un-intuitive. It seems impossible to extend the
involved analytic method to the branching
process with immigration. On the other hand, M\"orters and Ortgiese
\cite{MO08} provided a very useful probabilistic approach for
Theorem \ref{THM}, called the ``branching tree heuristic''
method. Our approach is built on the top of their powerful
arguments, and overcomes additional difficulties of immigration
effects. More specifically, we start with a fundamental decomposition for
$\cW$ given in (\ref{identity}). Then a suitable truncation is used in
order to handle the infinite series. To estimate the lower bound of
$\P(\cW\leq\varepsilon)$, we investigate when the least population
size happens. For the upper bound, we use the exponential
Chebyshev's inequality and estimate the Laplace transform of $\cW$.
The property of $\P(\cW\leq\varepsilon)$ is then obtained through
Tauberian type theorems.

Now we consider recursive distribution identities for
$(\cZ_n,n\geq0)$ satisfying $\cZ_0=Y_0$. For fixed integers $r\geq
0$ and $l \ge0$, let
$\xi_r(1),\ldots,\xi_r(\cZ_r)$ be the individuals in generation $r$,
and $\eta_{l}(j), j=1,\ldots,Y_l$ be the individuals of immigration
in generation $l$. Then for any $r \ge0$ and $n \ge r+1$,
\[
\cZ_n=\sum_{i=1}^{\cZ_r}Z_{n-r}
\bigl(\xi_r(i)\bigr)+\sum_{l=r+1}^{n}
\sum_{j=1}^{Y_l}Z_{n-l}\bigl(
\eta_{l}(j)\bigr).
\]
Here $(Z_n(v),n\geq0)$ is a supercritical G-W branching process
initiated with one individual $v$ and $W(v)$ is the limit of the
positive martingale $m^{-n}Z_n(v)$.

Dividing both sides of the above equality by $m^n$, then letting $n\to
\infty$, we get
%
\begin{equation}
\label{decomp} \cW=m^{-r}\sum_{i=1}^{\cZ_r}W
\bigl(\xi_r(i)\bigr)+ \sum_{l=r+1}^{\infty}m^{-l}
\sum_{j=1}^{Y_l}W\bigl(\eta_{l}(j)
\bigr).
\end{equation}
For simplicity, we rewrite
(\ref{decomp}) as
%
\begin{equation}
\label{identity} \cW=m^{-r}\sum_{i=1}^{\cZ_r}W_i+
\sum_{l=r+1}^{\infty}m^{-l}\sum
_{j=1}^{Y_l}W_l^j.
\end{equation}
Here all the
$W_i,W_l^j,i=1,\ldots,\cZ_r,l=r+1,\ldots,n,j=1,\ldots,Y_l$ are
independent and identically distributed as $W$. The relation (\ref
{identity}) is
the fundamental distribution identity of $\cW$ and it is used
repeatedly in our approach.

Next, we turn to consider a slightly different type of supercritical
branching process with immigration, which is denoted by $({\widetilde
\cZ}_n,n\geq0)$. The
only\vspace*{1pt} difference is to assume $\widetilde{\cZ}_0=1$. The corresponding
limit of $\widetilde{\cZ}_n/m^n$ is denoted by ${\widetilde\cW}$.
Then by
simple computation we get that
%
\begin{equation}
\label{decomp2}
{\widetilde\cW} =^d W + \frac{\cW}{m}
\end{equation}
in distribution, denoted by $=^d$ throughout this paper,
where $W$ and $\cW$ are independent.
Then owing to (\ref{decomp2})
and the fact that
%
\begin{eqnarray}
\P(W+\cW/m \leq\varepsilon)&\geq& \P(W \leq\varepsilon/2)\cdot\P(\cW
/m\leq
\varepsilon/2),
\nonumber\\[-8pt]\\[-8pt]
\P(W+\cW/m \leq\varepsilon)&\leq& \P(W\leq\varepsilon)\cdot\P(\cW/m
\leq
\varepsilon),\nonumber
\end{eqnarray}
we obtain the
following result as a consequence of combining Theorems \ref{THM} and
\ref{thm1}.
%
\begin{theorem}\label{thm2}
Assume the condition (\ref{XYlog}) holds.

\textup{(a)} If $p_0=0$, $p_1>0$ and $q_0>0$, then
\[
\P(\widetilde{\cW}\leq\varepsilon)\asymp\varepsilon^{|{\log
(p_1q_0)}|/\log m} \qquad\mbox{as }
\varepsilon\to0^+.
\]

\textup{(b)} If $p_0=0$, $p_1>0$ and $q_0=0$, then
\[
\log\P(\widetilde{\cW}\leq\varepsilon)\sim-\frac{K |{\log
p_1}|}{2(\log m)^2}|{\log
\varepsilon}|^2 \qquad\mbox{as } \varepsilon\to0^+
\]
with $K$ being defined as in Theorem \ref{thm1}\textup{(b)}.

\textup{(c)} If $p_0=0$ and $p_1=0$, then
\[
\log\P(\widetilde{\cW} \leq\varepsilon)\asymp-\varepsilon^{-\beta
/(1-\beta)}
\qquad\mbox{as } \varepsilon\to0^+
\]
with $\beta$ being defined as in Theorem \ref{THM}\textup{(b)}.

\textup{(d)} If $p_0>0$, then
\[
\P(\widetilde{\cW} \leq\varepsilon) \asymp\varepsilon^{|{\log
h(\rho)}|/\log m} \qquad\mbox{as }
\varepsilon\to0^+.
\]
\end{theorem}

Note that when $q_0=1$, that is, without immigration, Theorem \ref{thm2} recovers
Theorem \ref{THM}.

\section{\texorpdfstring{Proof of Theorem \protect\ref{thm1}: Lower bound}
{Proof of Theorem 2: Lower bound}}\label{lower}

We start with a simple but crucial probability estimate that is a
consequence of the condition $\E\log^+ Y <\infty$ in
(\ref{XYlog}).
%
\begin{lem}\label{lem} Under the condition that $\E\log^+ Y <\infty
$ in (\ref{XYlog}), for any fixed constant $\delta>0$, there exists
an integer $l$ such that
%
\begin{equation}
\label{lo2} \P\Bigl(\max_{i \ge l+1} Y_i \RMe ^{-\delta i}\le1
\Bigr) \ge \RMe^{-1}.
\end{equation}
\end{lem}
\begin{pf}
For any given $\delta>0$,
\begin{eqnarray*}
\sum_{i=1}^\infty\P\bigl(\log^+ Y \ge\delta i
\bigr) &=&\sum_{i=1}^\infty\sum
_{k= i}^\infty\P\bigl(k \le\delta^{-1} \log^+
Y < k+1\bigr)
\\
&=& \sum_{k=1}^\infty k \E\I\bigl(k \le
\delta^{-1} \log^+ Y < k+1\bigr)
\\
&\le& \delta^{-1} \E\log^+Y <\infty.
\end{eqnarray*}
Let $Y_i$ and $Y$ be our
independent and identically distributed immigration random
variables. Then for any large integer $l$ such that
%
\begin{equation}
\label{logYl} \sum_{i=l+1}^\infty\P\bigl(\log^+
Y \ge\delta i\bigr) \le1/2
\end{equation}
we have
\begin{eqnarray*}
\P\Bigl(\max_{i \ge l+1} Y_i\mathrm{e}^{-\delta i}\leq1\Bigr) &\ge&
\prod_{i=l+1}^\infty\bigl( 1-\P\bigl(\log^+ Y
\ge\delta i\bigr) \bigr)
\\
&\ge& \exp\Biggl( -2 \sum_{i=l+1}^\infty\P
\bigl(\log^+ Y \ge\delta i\bigr) \Biggr)
\\
&\ge& \RMe^{-1},
\end{eqnarray*}
here we used the fact that $(1-x)\RMe^{2x}$ is
increasing for $0 \le x < 1/2$. This finishes our proof of the lemma.
\end{pf}

\begin{pf*}{Proof of (a) and (b)}
For any $\varepsilon> 0$, let $k=k_\varepsilon$ be the integer such that
%
\begin{equation}
\label{defk} m^{-k}\leq\varepsilon< m^{-k+1},
\end{equation}
which is equivalent to saying
%
\begin{equation}
\label{defk1} k-1<|{\log\varepsilon}|/\log m\leq k \quad\mbox{or}\quad k=\bigl\lceil
|{\log
\varepsilon}|/\log m\bigr\rceil.
\end{equation}
Using the fundamental distribution identity (\ref{identity}) with
$r=0$, we have for a fixed integer $l$ to be chosen later,
%
\begin{eqnarray}
\label{est1} \P(\cW\leq\varepsilon)&=&\P\Biggl(\sum
_{i=0}^{\infty}m^{-i}\sum
_{j=1}^{Y_i}W_i^{j} \leq
\varepsilon\Biggr){}
\nonumber\\[-8pt]\\[-8pt]
&\geq&\P\Biggl(\sum_{i=0}^{k+l}m^{-i}
\sum_{j=1}^{Y_i}W_i^{j}
\leq\frac{\varepsilon}{2} \Biggr) \cdot\P\Biggl(\sum_{i=k+l
+1}^{\infty}m^{-i}
\sum_{j=1}^{Y_i}W_i^{j}
\leq\frac{\varepsilon}{2} \Biggr).\nonumber
\end{eqnarray}
The second term in (\ref{est1}) can be estimated by using $\varepsilon
\geq m^{-k}$ in (\ref{defk}) as below
%
\begin{eqnarray}
\label{llo} \P\Biggl(\sum_{i=k+l+1}^{\infty}m^{-i}
\sum_{j=1}^{Y_i}W_i^{j}
\leq\frac{\varepsilon}{2} \Biggr) &\geq& \P\Biggl(\sum
_{i=k+l+1}^{\infty}m^{-i}\sum
_{j=1}^{Y_i}W_i^{j} \leq
\frac{m^{-k}}{2} \Biggr)
\nonumber
\\
&=& \P\Biggl(\sum_{i=l+1}^{\infty}m^{-i}
\sum_{j=1}^{Y_i}W_i^{j}
\leq{1 \over2} \Biggr).
\end{eqnarray}
Note that the last equality follows from the independence and identical
distribution of all $W_i^{j}$'s and $Y_i$'s.
Next, we have by controlling the size of $Y_i$, $ i \ge l+1$, given
in Lemma~\ref{lem},
%
\begin{eqnarray}
\label{lo} &&\P\Biggl(\sum_{i=l+1}^{\infty}m^{-i}
\sum_{j=1}^{Y_i}W_i^{j}
\leq{1 \over2} \Biggr)
\nonumber
\\
&&\quad\geq \P\Biggl(\sum_{i=l
+1}^{\infty}m^{-i}
\sum_{j=1}^{Y_i}W_i^{j}
\leq{1 \over2}, \max_{i
\ge l+1} Y_i\mathrm{e}^{-\delta i}
\le1 \Biggr)
\\
&&\quad\geq\P\Biggl(\sum_{i=l +1}^{\infty}m^{-i}
\sum_{j=1}^{\lceil\exp
(\delta i)\rceil} W_i^{j}
\leq{1 \over2} \Biggr) \cdot\P\Bigl( \max_{i \ge l+1}
Y_i\mathrm{e}^{-\delta i} \le1 \Bigr).\nonumber
\end{eqnarray}
Using Chebyshev's inequality
for the first part of (\ref{lo}),
we get
%
\begin{eqnarray}
\P\Biggl(\sum_{i=l +1}^{\infty}m^{-i}
\sum_{j=1}^{\lceil\exp
(\delta i)\rceil} W_i^{j}
\leq{1 \over2} \Biggr) &\geq&1-2\E\sum_{i=l +1}^{\infty}m^{-i}
\sum_{j=1}^{\lceil\exp
(\delta i)\rceil} W_i^{j}
\nonumber\\[-8pt]\\[-8pt]
&\geq&1-2\sum_{i=l +1}^{\infty}m^{-i}
\bigl(\RMe^{\delta i}+1\bigr).\nonumber
\end{eqnarray}
For $\delta$ satisfying $\RMe^\delta<m$, we have $\sum_{i=l +1}^{\infty
}m^{-i}(\RMe^{\delta i}+1)<\infty$. Then we choose $\delta$ small enough
and integer $l$ large enough so that
%
\begin{equation}
\label{lo1} 2\sum_{i=l +1}^{\infty}m^{-i}
\bigl(\RMe^{\delta i}+1\bigr)<{1 \over2}.
\end{equation}
Combining (\ref{llo})--(\ref{lo1}) and Lemma \ref{lem}, we obtain
%
\begin{equation}
\label{sec} \P\Biggl(\sum_{i=k+l +1}^{\infty}m^{-i}
\sum_{j=1}^{Y_i}W_i^{j}
\leq\frac{\varepsilon}{2} \Biggr) \geq\P\Biggl(\sum_{i=l+1}^{\infty}m^{-i}
\sum_{j=1}^{Y_i}W_i^{j}
\leq{1 \over2} \Biggr)\geq{1 \over2\RMe}.
\end{equation}

Now back to the first part of (\ref{est1}), we have to handle it
under conditions (a) and (b) separately. In the case (a) with $q_0>0$, we
have the simple estimate
%
\begin{equation}
\label{estimation2} \P\Biggl(\sum_{i=0}^{k+l}m^{-i}
\sum_{j=1}^{Y_i}W_i^{j}
\leq\frac{\varepsilon}{2} \Biggr) \geq\P(Y_0=\cdots=Y_{k+l}=0
)=q_0^{k+l+1}.
\end{equation}
Using $k-1< |{\log\varepsilon}|/\log m$ in (\ref{defk1}), it's easy to
deduce that
%
\begin{equation}
\label{estimation13} q_0^{k}\geq q_0 \cdot
q_0^{|{\log\varepsilon}|/\log m}=q_0\varepsilon^{|{\log
q_0}|/\log m}.
\end{equation}
Combining (\ref{est1}) and
(\ref{sec})--(\ref{estimation13}), we have shown the lower bound in
Theorem \ref{thm1}(a).

For the case (b) with $q_0=0$, we have, recalling the definition of
$K=\inf\{ n\dvtx  q_n >0\}$,
%
\begin{eqnarray}
\label{up3}
&&
\P\Biggl(\sum_{i=0}^{k+l}m^{-i}
\sum_{j=1}^{Y_i}W_{i}^j
\leq\frac
{\varepsilon}{2} \Biggr)\nonumber\\
&&\quad\geq\P\Biggl(\sum_{i=0}^{k+l}m^{-i}
\sum_{j=1}^{Y_i}W_{i}^j
\leq\frac{\varepsilon}{2},Y_0=\cdots=Y_{k+l}=K \Biggr)
\\
&&\quad=\P\Biggl(\sum_{i=0}^{k+l}m^{-i}
\sum_{j=1}^{K}W_{i}^j
\leq\frac{\varepsilon}{2} \Biggr) \cdot q_K^{k+l+1}.\nonumber
\end{eqnarray}
The above probability of sums can be bounded termwise, and thus
%
\begin{eqnarray}
\label{up31}
&&\P\Biggl(\sum_{i=0}^{k+l}m^{-i}
\sum_{j=1}^{K}W_{i}^j
\leq\frac{\varepsilon}{2} \Biggr)\nonumber\\
&&\quad\geq\P\biggl( \max_{0 \le i \le
k+l}
\max_{1 \le j \le K} m^{-i}W_i^j\leq
\frac{\varepsilon/2}{K(
k+l+1)} \biggr)
\nonumber\\[-8pt]\\[-8pt]
&&\quad=\prod_{i=0}^{k+l}\P^K
\biggl(m^{-i}W\leq\frac{\varepsilon
/2}{K(k+l+1)} \biggr)
\nonumber\\
&&\quad\ge \prod_{i=0}^{k+l}\P^K
\biggl(W\leq\frac{m^{i-k}/2}{K (k+
l+1)} \biggr),\nonumber
\end{eqnarray}
%
where we use the independent and identically distributed property of
all $W_i^j$'s in the last equality and $\varepsilon\geq m^{-k}$ from
(\ref{defk}) in the last inequality.

From Theorem \ref{THM}(a) there exists a constant $c>0$
such that, for $i=0,1,\ldots, k+l$,
%
\begin{equation}
\label{up16} \P\biggl(W\leq\frac{m^{i-k}/2}{K(k+l+1)} \biggr)\geq c
\biggl(
\frac{m^{i-k}/2}{K(k+l+1)} \biggr)^{|{\log p_1}|/\log m}.
\end{equation}
Combining (\ref{est1}), (\ref{sec}) and (\ref{up3})--(\ref{up16})
together, and taking summation over $0 \le i \le k+l$ after
taking logarithm, we have
\begin{eqnarray*}
\log\P(\cW\leq\varepsilon) &\geq&-\frac{K |{\log p_1}|}{2}k^2-\RMO(k\log k)
\\
&\ge& -\frac{K |{\log p_1}|}{2(\log m)^2}|{\log\varepsilon}|^2-\RMO\bigl(\log
\varepsilon^{-1}\log\log\varepsilon^{-1}\bigr),
\end{eqnarray*}
which follows easily from $k<1+ |{\log\varepsilon}|/\log m$ in (\ref{defk1}).
\noqed\end{pf*}

\begin{pf*}{Proof of (c)} First observe that, in this setting
with $\gamma=\inf\{n\dvtx  p_n >0\}\geq2, K=\inf\{n\dvtx  q_n >0\}\geq1$,
the smallest number of particles in generation $n$ $(n\geq1)$ is
%
\begin{equation}
\label{bn} b(n):=K\bigl(\gamma^n+\gamma^{n-1}+\cdots+1
\bigr) =K\bigl(\gamma^{n+1}-1\bigr)/(\gamma-1).
\end{equation}
It is also easy to see that the chance this occurs is
%
\begin{equation}
\label{prob1} \P\bigl(\cZ_n=b(n)\bigr) =p_\gamma^{b(n-1)+\cdots+b(0)}
q_K^{n+1}:=p_\gamma^{B(n)}
q_K^{n+1},
\end{equation}
where
%
\begin{equation}
\label{defB} B(0)=0,\qquad B(n)=b(n-1)+\cdots+b(0)=\frac{K (\gamma
^{n+1}-(n+1)\gamma
+n )}{(\gamma-1)^2},\qquad n\geq1.\quad
\end{equation}
Given $\varepsilon> 0$, we can choose $k=k_\varepsilon$ such that
%
\begin{equation}
\label{ineq} \frac{\gamma^{k}}{m^{k}}\leq\varepsilon< \frac{
\gamma^{k-1}}{m^{k-1}},
\end{equation}
which is equivalent to saying
%
\begin{equation}
\label{ineq1} k-1 < |{\log\varepsilon}|/ \log(m/ \gamma) \leq k\quad \mbox
{or}\quad k= \bigl
\lceil|{\log\varepsilon}|/ \log(m/ \gamma) \bigr\rceil.
\end{equation}
Next, let $l$ be an integer that will be determined later.
Using the fundamental distribution identity (\ref{identity}) with
$r=k+l$ and (\ref{prob1}), we have
%
\begin{eqnarray}
\label{upp1} &&\P(\cW\leq\varepsilon)
\nonumber
\\
&&\quad\geq\P\bigl(\cW\leq(\gamma/ m)^{k}|\cZ_{k+l}=b(k+l)
\bigr) \P\bigl(\cZ_{k+l}=b(k+l)\bigr){}
\nonumber\\[-8pt]\\[-8pt]
&&\quad=\P\Biggl(m^{-k-l}\sum_{i=1}^{b(k+l)}W_i+
\sum_{i=k+l+1}^{\infty}m^{-i}\sum
_{j=1}^{Y_i}W_i^j
\leq(\gamma/m)^{k} \Biggr) p_\gamma^{B(k+l)}
q_K^{k+l+1}
\nonumber
\\
&&\quad\ge \P\Biggl(\sum_{i=1}^{b(k+l)}W_i
\leq{m^l\gamma^{k} \over2} \Biggr)\P\Biggl(\sum_{i=1}^{\infty}m^{-i}
\sum_{j=1}^{{Y_i}}W_i^j
\leq{m^l \gamma^{k} \over2} \Biggr) p_\gamma^{B(k+l)}
q_K^{k+l+1}.\nonumber
\end{eqnarray}
For the first term in (\ref{upp1}) we have by
Chebyshev's inequality
and choosing suitable $l$
%
\begin{eqnarray}
\label{upp3} \P\Biggl(\sum_{i=1}^{b(k+l)}W_i
\leq m^l\gamma^{k}/2 \Biggr) &\ge& 1-{2 \over m^l\gamma^k}
\E\sum_{i=1}^{b(k+l)} W_i \nonumber\\
&=& 1-
{2
b(k+l) \over m^l\gamma^k}
\\
&\ge& 1-{2K \gamma\over\gamma-1} (\gamma/m)^l \ge1/2,\nonumber
\end{eqnarray}
where $\E W=1$ and $b(n) \le K(\gamma-1)^{-1}
\gamma^{n+1}$ from (\ref{bn}) are used.

For the second part of (\ref{upp1}), we have
%
\begin{eqnarray}
\label{upp5} \P\Biggl(\sum_{i=1}^{\infty}m^{-i}
\sum_{j=1}^{{Y_i}}W_i^j
\leq\frac{m^l
\gamma^{k}}{2} \Biggr) &=&\P\Biggl(\sum_{i=l+1}^{\infty}m^{-i}
\sum_{j=1}^{{Y_i}}W_i^j
\leq\frac{
\gamma^{k}}{2} \Biggr)
\nonumber\\
&\geq&\P\Biggl(\sum_{i=l+1}^{\infty}m^{-i}
\sum_{j=1}^{{Y_i}}W_i^j
\leq\frac{
1}{2} \Biggr)\\
&\geq& \RMe^{-1}/2,\nonumber
\end{eqnarray}
where the last inequality follows from (\ref{sec}).

Combing (\ref{upp1})--(\ref{upp5}), we get
%
\begin{equation}
\label{in2} \P(\cW\leq\varepsilon)\geq p_\gamma^{B(k+l)}
q_K^{k+l+1}\RMe^{-1}/4.
\end{equation}
Recalling the definition of $B(k+l)$ in (\ref{defB})
and $k-1 < |{\log\varepsilon}|/\log(m/\gamma)$ in (\ref{ineq1}), we see
\[
B(k+l)\leq\frac{K}{(\gamma-1)^2}\gamma^{k+l+1}\leq C \gamma^{|{\log
\varepsilon}|/\log(m/\gamma)}=C
\varepsilon^{-\beta/(1-\beta)},
\]
where $\beta$ is defined as in Theorem \ref{THM}(b) and $C$ is a
positive
constant. Therefore from (\ref{in2}), we obtain
\[
\log\P(\cW\leq\varepsilon)\geq-C \varepsilon^{-\beta/(1-\beta)}
\]
for some constant $C>0$.
\noqed\end{pf*}

\section{\texorpdfstring{Proof of Theorem \protect\ref{thm1}: Upper bound}
{Proof of Theorem 2: Upper bound}}\label{upper}
As we can see from the arguments in Section \ref{lower}, only the
finite terms in (\ref{identity}) are contributing to the small value
probabilities of $\cW$. Hence, we take only $r=0$ in
(\ref{identity}), choose a suitable cut off~$k$, and focus on
properties of $\sum_{l=0}^{k}m^{-l}\sum_{j=1}^{Y_l}W_l^j$.

\begin{pf*}{Proof of (a)} Let $k=k_\varepsilon$ be the
integer defined as in (\ref{defk}). Using the fundamental
distribution identity (\ref{identity}) with $r=0$ and exponential
Chebyshev's inequality,
we have
%
\begin{eqnarray}
\label{upper1} \P(\cW\leq\varepsilon) &\leq&\P\Biggl(\sum
_{i=0}^{k}m^{-i}\sum
_{j=1}^{Y_i}W_i^{j} \leq
\varepsilon\Biggr){}
\nonumber\\[-8pt]\\[-8pt]
&\leq&\RMe^{\lambda\varepsilon}\cdot\E\exp\Biggl(-\lambda\sum
_{i=0}^{k}m^{-i}\sum
_{j=1}^{Y_i}W_i^{j} \Biggr)
\qquad\mbox{for any } \lambda>0.\nonumber
\end{eqnarray}
Noticing that all the
$(W_i^{j}, i=0,\ldots,k,j=1,\ldots,Y_i)$ are independent, we have
%
\begin{equation}
\label{upper2} \E\exp\Biggl(-\lambda\sum_{i=0}^{k}m^{-i}
\sum_{j=1}^{Y_i}W_i^{j}
\Biggr)=\prod_{i=0}^{k}\E\exp\Biggl(-
\lambda m^{-i}\sum_{j=1}^{Y_i}W_i^{j}
\Biggr).
\end{equation}
Conditioning on $Y_i=0$ or $Y_i\geq1$, we have
%
\begin{equation}
\label{upper11} \E\exp\Biggl(-\lambda m^{-i}\sum
_{j=1}^{Y_i}W_i^{j} \Biggr)\le
q_0+(1-q_0)\E\exp\bigl(-\lambda m^{-i}
W_i^{1} \bigr)\leq q_0(1+
\delta_i),
\end{equation}
where
%
\begin{equation}
\label{bound11} \delta_i=q_0^{-1}\E\exp
\bigl(-\lambda m^{-i} W_i^{1}
\bigr)=q_0^{-1}\E\exp\bigl(-\lambda m^{-i} W
\bigr),\qquad i=0,\ldots,k.
\end{equation}
Substituting (\ref{upper11}) into (\ref{upper1}) and
letting $\lambda=\varepsilon^{-1}$, we obtain
\[
\P(\cW\leq\varepsilon)\leq \RMe q_0^{k+1}\prod
_{i=0}^{k}(1+\delta_i).
\]
Since $k\geq|{\log\varepsilon}|/\log m$ in (\ref{defk1}), we have
\[
q_0^k \le\varepsilon^{|{\log q_0}|/\log m}.
\]
So we finish the proof by showing
%
\begin{equation}
\label{bound} \sum_{i=0}^{k}\log(1+
\delta_i)\leq\sum_{i=0}^{k}
\delta_i \leq M,
\end{equation}
where $M>0$ is a constant independent of
$\varepsilon$ (noticing that the $k$ depends on $\varepsilon$).
\noqed\end{pf*}

In order to estimate $\delta_i$, we need the following fact given in
Li \cite{L12+}.
%
\begin{lem}\label{Taub}
\textup{(i)} Assume $V$ is a positive random variable and $\alpha>0$ is a
constant. Then
\[
\P(V \leq t)\leq C_1 t^{\alpha} \qquad\mbox{for some constant
$C_1>0$ and all $t >0$}
\]
is equivalent to
\[
\E \RMe^{-\lambda V}\leq C_2 \lambda^{-\alpha} \qquad\mbox{for some
constant $C_2>0$ and all $\lambda>0$}.
\]

\textup{(ii)} Assume V is a positive random variable and $\alpha>0,\theta
\in
\R$, or $\alpha=0,\theta>0$ are constants. Then we have
\[
\log\P(V \leq t)\leq-C_1t^{-\alpha}|{\log t}|^\theta
\qquad\mbox{for some constant $C_1>0$ and all $t >0$}
\]
is equivalent to
\[
\log\E \RMe^{-\lambda V}\leq-C_2\lambda^{\alpha/(1+\alpha)}(\log
\lambda)^{\theta/(1+\alpha)} \qquad\mbox{for some constant $C_2>0$ and all $
\lambda>0$}.
\]
\end{lem}

To show (\ref{bound}), we have to argue
separately according to $p_1>0$ or $p_1=0$.
When $p_1>0$, by Theorem \ref{THM}(a) and Lemma \ref{Taub}(i),
there exists a constant $C>0$ satisfying that
%
\begin{equation}
\label{exp-v} \E \RMe^{-\lambda W}\leq C \lambda^{-|{\log p_1}|/\log
m},\qquad
\lambda>0.
\end{equation}
Combining (\ref{bound11}) with $\lambda=\varepsilon^{-1}$, then
using (\ref{exp-v}), we have
\begin{eqnarray*}
\sum_{i=0}^{k}\delta_i &=&
q_0^{-1} \sum_{i=0}^k
\E\exp\bigl( -\varepsilon^{-1} m^{-i} W \bigr)
\\
&\leq& q_0^{-1} C \sum_{i=0}^k
\bigl(\varepsilon m^i\bigr)^{|{\log p_1}|/\log m}
\\
&=& Cq_0^{-1} \varepsilon^{|{\log p_1}|/\log m } \sum
_{i=0}^kp_1^{-i}
\\
&\le& C' \varepsilon^{|{\log p_1}|/\log m} \cdot p_1^{-k}
\leq C'p_1^{-1},
\end{eqnarray*}
where $C'$ is a constant and the last inequality follows from (\ref{defk1}).

When $p_1=0$, using Theorem \ref{THM}(b) and Lemma \ref{Taub}(ii)
with $\alpha=\beta/(1-\beta)$ and $\theta=0$, we have for
some constant $b>0$,
%
\begin{equation}
\label{bound2} \log\E \RMe^{-\lambda W}\leq-b\lambda^{\beta}, \qquad\lambda>0,
\end{equation}
from which it's similar to show that (\ref{bound}) holds. Indeed,
setting $\lambda=\varepsilon^{-1}$ in (\ref{bound11}), and then
using (\ref{bound2}) and $\varepsilon< m^{-k+1}$ from (\ref{defk}),
we obtain
\begin{eqnarray*}
\sum_{i=0}^{k}\delta_i &=&
q_0^{-1} \sum_{i=0}^k
\E\exp\bigl( -\varepsilon^{-1} m^{-i} W \bigr)
\\
&\leq& q_0^{-1} \sum_{i=0}^k
\exp\bigl(-b\varepsilon^{-\beta} m^{-i\beta
}\bigr)
\\
&\leq& q_0^{-1} \sum_{i=0}^k
\exp\bigl(-bm^{(k-i-1)\beta}\bigr)
\\
&\leq&q_0^{-1} \sum_{i=0}^\infty
\exp\bigl(-bm^{(i-1)\beta}\bigr)<\infty.
\end{eqnarray*}

\begin{pf*}{Proof of (b)}
Let $k$ be defined as in
(\ref{defk}). Using (\ref{upper1}) and $Y_i\geq K$ for
any $i\geq0$,
%
\begin{equation}
\label{low1} \P(\cW\leq\varepsilon)\leq \RMe^{\lambda\varepsilon}\prod
_{i=0}^{k}\prod_{j=1}^K
\E\exp\bigl(-\lambda m^{-i}W_{i}^j \bigr),\qquad
\lambda>0.
\end{equation}
In the case (b) with $p_1>0$, substituting (\ref{exp-v}) into
(\ref{low1}) with $\lambda=\varepsilon^{-1}$, we obtain
\[
\P(\cW\leq\varepsilon)\leq \RMe\prod_{i=0}^{k}
\prod_{j=1}^K C\bigl(\varepsilon
m^{i}\bigr)^{|{\log p_1}|/\log m}.
\]
Taking the logarithm we obtain
\begin{eqnarray*}
\log\P(\cW\leq\varepsilon)&\leq& 1+K (k+1) \bigl(\log C-|{\log\varepsilon}
|\cdot|{\log
p_1}|/\log m\bigr)+k(k+1)\cdot K|{\log p_1}|/2
\\
&=&-k\cdot|{\log\varepsilon}|\cdot K|{\log p_1}|/\log m +
(k-1)^2\cdot K|{\log p_1}|/2+ \RMO(k)
\\
&\leq&-\frac{K |{\log p_1}| }{2(\log m)^2}|{\log\varepsilon}|^2+\RMO\bigl(|{\log
\varepsilon}|\bigr),
\end{eqnarray*}
where the last inequality follows from $k-1<|{\log\varepsilon}|/\log
m\leq k$, which is given in (\ref{defk1}).
\noqed\end{pf*}

\begin{pf*}{Proof of (c)} It is clear that
%
\begin{equation}
\P(\cW\leq\varepsilon)\leq\P(W\leq\varepsilon),
\end{equation}
and therefore we finish the proof of $(c)$ by using estimate in Theorem
\ref{THM}(b).
\noqed\end{pf*}

\section{\texorpdfstring{Proof of Theorem \protect\ref{thm1}(d)}
{Proof of Theorem 2(d)}}\label{p0}

If $p_0>0$, then $f(s)=s$ has a unique solution $\rho\in(0,1)$ and
$\P(W=0)=\rho$. By means of the Harris--Sevastyanov
transformation
\[
\widetilde{f}(s):=\frac{f((1-\rho)s+\rho)-\rho}{(1-\rho)},
\]
$\widetilde f$ defines\vspace*{1pt} a new branching mechanism with $\widetilde
{p}_0=0$ and
$\widetilde{f}'(1)=m$.
We use $(\widetilde{Z}_n,n\geq0)$ to denote the
corresponding branching process and $\widetilde{W}$
to denote the limit of
$m^{-n}\widetilde{Z}_n$. By Theorem 3.2 in~\cite{H48},
%
\begin{equation}
\label{hstransform}W=^{d}W_0\cdot\widetilde{W},
\end{equation}
where $W_0$ is
independent of $\widetilde{W}$ and takes the values 0 and
$1/(1-\rho)$ with probabilities $\rho$ and $1-\rho$, respectively.
Notice that the small value probability of $\widetilde{W}$ has the
asymptotic behavior described in Theorem \ref{THM}(a) with
$\widetilde{p}_1=\widetilde{f}'(0)=f'(\rho)>0$, and $\tau=|{\log
\widetilde{p}_1}|/\log m$, that is,
%
\begin{equation}
\label{hssmall} \P( \widetilde{W}\leq\varepsilon) \asymp\varepsilon
^{\tau}.
\end{equation}

Now we start to prove the lower bound.
For any $\varepsilon> 0$, let $k=k_\varepsilon$ be the integer
defined in~(\ref{defk}). Then using (\ref{est1}) and (\ref{sec}), we
only need
to estimate the first part of (\ref{est1}):
%
\begin{eqnarray}
\label{esti} \P\Biggl(\sum_{i=0}^{k+l}m^{-i}
\sum_{j=1}^{Y_i}W_i^{j}
\leq\frac{\varepsilon}{2} \Biggr) &\geq&\prod_{i=0}^{k+l}
\P\Biggl(\sum_{j=1}^{Y_i}W_i^{j}=0
\Biggr)
\nonumber\\[-8pt]\\[-8pt]
&=&\prod_{i=0}^{k+l} \Biggl(\sum
_{n=0}^{\infty}q_n \P^n(W=0)
\Biggr)=h(\rho)^{k+l+1},\nonumber
\end{eqnarray}
where $h$ is the generating function of immigration $Y$.
Using $k-1< |{\log\varepsilon}|/\log m$
given in (\ref{defk1}), it's easy
to deduce that
%
\begin{equation}
\label{est13} h(\rho)^{k}\geq h(\rho)\cdot h(\rho)^{|{\log\varepsilon
}|/\log
m}=h(
\rho)\cdot\varepsilon^{|{\log h(\rho)}|/\log m}.
\end{equation}
Combining (\ref{est1}), (\ref{sec}), (\ref{esti}) and (\ref{est13}),
we obtain the lower bound of (d).

Next, we show the upper bound. Using (\ref{hstransform}),
we have
%
\begin{equation}
\E \RMe^{-\lambda W}=\rho+\E \RMe^{-\lambda
W}\mathbb{I}_{\{W>0\}}:=\rho+
\delta(\lambda),\qquad \lambda>0.
\end{equation}
%
Using (\ref{upper1}), (\ref{upper2}) and the independent and
identically distributed property of all the
$(W_i^{j}, i=0,\ldots,k,j=1,\ldots,Y_i)$, we have
%
\begin{eqnarray}
\label{uppp2} \P(\cW\leq\varepsilon) &\leq&\RMe^{\lambda\varepsilon}\prod
_{i=0}^{k}h \bigl(\rho+\delta\bigl(\lambda
m^{-i}\bigr) \bigr)
\nonumber\\[-8pt]\\[-8pt]
&=&\bigl(h(\rho)\bigr)^{k+1}\exp\Biggl(\lambda\varepsilon+\sum
_{i=0}^{k}\log\bigl(h \bigl(\rho+\delta
\bigl(\lambda m^{-i}\bigr) \bigr)/h(\rho) \bigr) \Biggr),\nonumber
\end{eqnarray}
where $\lambda=\lambda_k$
depends on $k(=k_\varepsilon)$
and is given later.
Since $k\geq|{\log\varepsilon}|/\log m$
from (\ref{defk1}), we have
%
\begin{equation}
\label{uppp3} \bigl(h(\rho)\bigr)^k \le\varepsilon^{|{\log h(\rho)}|/\log m}.
\end{equation}
Next we show there is a constant $M>0$, which does not depend on
$\varepsilon$, such that
%
\begin{eqnarray}
\label{finite}&&\lambda\varepsilon+\sum_{i=0}^{k}
\log\bigl(h \bigl(\rho+\delta\bigl(\lambda m^{-i}\bigr) \bigr)/h(\rho)
\bigr)
\nonumber\\[-8pt]\\[-8pt]
&&\quad\leq \lambda m^{-k+1}+h(\rho)^{-1}\sum
_{i=0}^{k} \bigl(h \bigl(\rho+\delta\bigl(\lambda
m^{-i}\bigr) \bigr)-h(\rho) \bigr)\le M.\nonumber
\end{eqnarray}
Since $\delta(\lambda m^{-x})$ is increasing with respect to $x$, we have
%
\begin{equation}
\label{sum2a}
\sum_{i=0}^{k} \bigl(h \bigl(
\rho+\delta\bigl(\lambda m^{-i}\bigr) \bigr)-h(\rho) \bigr)\leq\int
_{0}^{k+1} \bigl(h \bigl(\rho+\delta\bigl(\lambda
m^{-x}\bigr) \bigr)-h(\rho) \bigr)\,\mathrm{d}x.
\end{equation}
Note that $\delta(\lambda)=(1-\rho)\E \RMe^{-(\lambda/(1-\rho))
\widetilde W}$. By (\ref{hssmall}) and Lemma \ref{Taub}(i), there
exists a constant $C>0$
such that
%
\begin{equation}
\delta\bigl(\lambda m^{-x}\bigr)\leq C \bigl(\lambda m^{-x}
\bigr)^{-\tau}
\end{equation}
with $\tau=|{\log f'(\rho)}|/\log m$. Thus, we have
%
\begin{eqnarray}
\label{sum2} &&\sum_{i=0}^{k} \bigl(h
\bigl(\rho+\delta\bigl(\lambda m^{-i}\bigr) \bigr)-h(\rho) \bigr)
\nonumber
\\
&&\quad\leq\int_{0}^{k+1} \bigl(h \bigl(\rho+C\bigl(
\lambda m^{-x}\bigr)^{-\tau
} \bigr)-h(\rho) \bigr)\,\mathrm{d}x
\nonumber\\[-8pt]\\[-8pt]
&&\quad=1/(\tau\log m)\cdot\int_{\lambda^{-\tau}}^{\lambda^{-\tau}
m^{(k+1)\tau}}1/y\cdot
\bigl(h (\rho+Cy )-h(\rho) \bigr)\,\mathrm{d}y
\nonumber
\\
&&\quad\le1/(\tau\log m)\cdot\int_{0}^{\lambda^{-\tau}
m^{(k+1)\tau}}1/y\cdot
\bigl(h (\rho+Cy )-h(\rho) \bigr)\,\mathrm{d}y.\nonumber
\end{eqnarray}
As $\rho<1$, we may choose $\delta_0>0$ such that $\rho+\delta_0<1$.
Next, we choose $\lambda= ({C}/{\delta_0} )^{1/\tau
} m^{(k+1)}$ in order to assure
$\rho+Cy<1$ so that $h (\rho+Cy )$ is well defined.
Indeed, we have
%
\begin{equation}
\label{sum1}\lambda m^{-k+1}=m^2 ({C}/{\delta_0}
)^{1/\tau}:=M_1
\end{equation}
and
\[
\rho+Cy\le\rho+C\lambda^{-\tau} m^{(k+1)\tau}= \rho+
\delta_0<1,\qquad y\le\lambda^{-\tau} m^{(k+1)\tau}.
\]
Then we follow (\ref{sum2}) to get
%
\begin{eqnarray}\label{sum2b}
\sum_{i=0}^{k} \bigl(h\bigl(\rho+\delta\bigl(\lambda m^{-i}\bigr) \bigr)-h(\rho) \bigr)
&\le&1/(\tau\log m)\cdot\int_{0}^{\delta_0/C}1/y\cdot
\bigl(h (\rho+Cy )-h(\rho) \bigr)\,\mathrm{d}y\nonumber\\[-8pt]\\[-8pt]
&:=&M_2<\infty,\nonumber
\end{eqnarray}
where we used
\[
\lim_{y\to
0}1/y\cdot\bigl(h (\rho+Cy )-h(\rho) \bigr)
=Ch'(\rho)<\infty.
\]
From (\ref{finite}), (\ref{sum1}) and (\ref{sum2b}),
we obtain that (\ref{finite}) holds with $M=M_1+M_2$, and finish the
proof of Theorem \ref{thm1}(d).

\section*{Acknowledgements}

We would like to thank referees for insightful remarks and suggestions.
The second author is supported in part by DMS-08-05929, DMS-11-06938,
NSFC-10928103. The third author is supported in part by NSFC (Grants
10971003 and 11128101) and Specialized Research Fund for the
Doctoral Program of Higher Education of China.



\printhistory

\end{document}